\documentclass{amsart}

\usepackage[utf8x]{inputenc}
\usepackage{amsfonts,amsmath,amssymb}
\usepackage{mathtools}

\DeclareMathOperator{\cl}{cl}
\DeclareMathOperator{\non}{non}
\DeclareMathOperator{\cov}{cov}
\DeclareMathOperator{\cof}{cof}
\DeclareMathOperator{\add}{add}
\DeclareMathOperator{\Con}{Con}

\newcommand{\concat}{
  \text{\raisebox{.8ex}{$\smallfrown$}}
}
\newtheorem{theorem}{Theorem}
\newtheorem{lemma}[theorem]{Lemma}

\newtheorem{example}[theorem]{Example}
\newtheorem{definition}[theorem]{Definition}

\newtheorem{fact}[theorem]{Fact}
\newtheorem{corollary}[theorem]{Corollary}

\begin{document}

\title{Topology idealized}

\author[G. We\l{}yczko]{Grzegorz We\l{}yczko}
\email{277511@student.pwr.edu.pl}

\author[Sz. \.Zeberski]{Szymon \.Zeberski}       \email{szymon.zeberski@pwr.edu.pl}

\address[G. We\l{}yczko, Sz. \.Zeberski]{Faculty of Pure and Applied Mathematics, Wrocław University of Science and Technology, Wybrzeże Stanisława Wyspiańskiego 27, 50-370 Wrocław, Poland}

\subjclass{03E017; 03E75, 26A03, 03E35.}
\keywords{local function, closure operator, cardinal invariant, Polish space, ideal, continuous function}

\begin{abstract}
Inspired by recent work of A. Njamcul and A. Pavlovi\'c, we study topological spaces expanded by ideals of sets. We investigate the relationship between the cardinal invariants of the ideal topological space $(X, \tau_\mathcal{I})$, the underlying topology $(X, \tau)$, and the cardinal characteristics of the ideal $\mathcal{I}$. Furthermore, we analyze continuous maps between standard Polish spaces and their ideal expansions. We establish rigidity phenomena showing that continuous functions from $(X, \tau)$ into $(X, \tau_\mathcal{I})$ are severely constrained (e.g., taking at most countably many values on Polish spaces and being strictly constant on $\mathbb{R}$), whereas continuous functions from $(\mathbb{R}, \tau_\mathcal{I})$ to $(\mathbb{R}, \tau)$ collapse to standard continuous functions under mild conditions on the ideal.
\end{abstract}

\maketitle

\section{Inspirations and basic definitions}

The study of topology expansion via ideals dates back to the classical work of Kuratowski \cite{Kuratowski} and Vaidyanathaswamy \cite{Vaid}, and was later systematically developed by Janković and Hamlett. Recently, a renewed interest in ideal-generated topologies was sparked by the work of A. Njamcul and A. Pavlovi\'c \cite{Cale}. Motivated by their results, we aim to explore the structural, cardinal, and function-space properties of such ideal expansions.

Throughout this paper, we adopt standard notation from general topology and set theory. For topological spaces $(X, \tau)$ and $(Y, \rho)$, the symbol $C((X, \tau), (Y, \rho))$ denotes the set of all continuous functions from $X$ to $Y$, that is:
\[
C((X, \tau), (Y, \rho)) = \{f \in Y^X : \forall U \in \rho \;\; f^{-1}[U] \in \tau\}.
\]
We will also write ``$f \colon (X, \tau) \to (Y, \rho)$ is continuous'' to state that $f \in C((X, \tau), (Y, \rho))$. For a point $x \in X$, the symbol $\tau(x)$ denotes the family of all neighborhoods of $x$. Recall that a family $\mathcal{I} \subseteq \mathcal{P}(X)$ is an ideal on $X$ if it is closed under taking subsets and under finite unions.

Let us recall the foundational definitions and concepts considered in \cite{Cale}.

\begin{definition}
Let $(X,\tau)$ be a topological space and let $\mathcal{I}$ be an ideal of subsets of $X$.
\begin{itemize}
    \item For any $A\subseteq X$, the \emph{local function} $A^*$ associated with $\mathcal{I}$ and $\tau$ is defined as
    \[
    A^* = \{x\in X : (\forall U\in\tau)(x\in U \rightarrow U\cap A\notin \mathcal{I})\}.
    \]
    \item The operator $\cl^*(A) = A \cup A^*$ defines a Kuratowski closure operator on $X$.
    \item The topology on $X$ generated by $\cl^*$ is denoted by $\tau_\mathcal{I}$.
\end{itemize}
\end{definition} 

\begin{fact}
A basis for the expanded topology $\tau_\mathcal{I}$ consists of all sets of the form $U\setminus I$, where $U\in\tau$ and $I\in\mathcal{I}$. 
\end{fact}

Throughout this paper, we focus primarily on ideal spaces of the form $(X,\tau_\mathcal{I})$, where $X$ is a Polish space and $\mathcal{I}$ is a non-trivial $\sigma$-ideal containing all countable subsets of $X$. 

We adopt standard notation for topological cardinal invariants: $w(X, \tau)$, $d(X, \tau)$, $\chi(X, \tau)$, $t(X, \tau)$, and $L(X, \tau)$ denote the weight, density, character, tightness, and Lindel\"of degree of $(X, \tau)$, respectively (see e.g., Engelking \cite{Engelking} or Juhász \cite{Juhasz}). For an ideal $\mathcal{I}$ on $X$, $\add(\mathcal{I})$, $\cov(\mathcal{I})$, $\non(\mathcal{I})$, and $\cof(\mathcal{I})$ represent the standard cardinal invariants of the ideal, i.e.
\begin{align*}
        \add(\mathcal{I})&=\min \{ |\mathcal{A}|:\; \mathcal{A}\subseteq \mathcal{I}\land \bigcup\mathcal{A}\notin \mathcal{I} \};
        \\
        \non(\mathcal{I})&=\min \{ |A|:\; A\subseteq X\land A\notin \mathcal{I}\};
        \\
        \cov(\mathcal{I})&=\min \{ |\mathcal{A}|:\; \mathcal{A}\subseteq \mathcal{I}\land \bigcup\mathcal{A} = X\};
        \\
        \cof(\mathcal{I})&=\min \{ |\mathcal{A}|:\; \mathcal{A}\subseteq \mathcal{I}\land (\forall A\in \mathcal{I})(\exists B\in\mathcal{A}) (A\subseteq B)\}.
    \end{align*}

We will use the following elementary fact regarding minimal bases.

\begin{lemma}\label{lem:base-weight}
Let \((X, \tau)\) be a topological space and let \(\mathcal{B}\) be a base for \((X, \tau)\). Then there exists a subfamily \(\mathcal{B}' \subseteq \mathcal{B}\) such that \(|\mathcal{B}'| = w(X, \tau)\) and \(\mathcal{B}'\) is still a base for \((X, \tau)\).
\end{lemma}

In Section 2, we examine cardinal invariants of the space $(X, \tau_\mathcal{I})$ in relation to those of $(X, \tau)$ and $\mathcal{I}$, establishing bounds on density, tightness, and weight, as well as several consistency results. In Section 3, we analyze continuous mappings between standard and ideal topologies, demonstrating rigidity results for continuous functions defined on Polish spaces.

\section{Cardinal invariants}

We start this section with the following easy observation.

\begin{theorem}
Let \((X, \tau)\) be a topological space and let \(\mathcal{I}\) be an ideal on \(X\). Then \(\non(\mathcal{I}) \leq d(X, \tau_\mathcal{I})\).
\end{theorem}
\begin{proof}
    Suppose that \(D\) is a set of cardinality strictly less than \(\non(\mathcal{I})\). Then \(D \in \mathcal{I}\), hence \(X \setminus D\) is an open set in \((X, \tau_\mathcal{I})\) that has an empty intersection with \(D\).
\end{proof}

The following example shows that in many natural cases the other inequality also holds.

\begin{example}
    Let \(\mathcal{I}\) be a \(\sigma\)-ideal on \(\mathbb{R}\) such that for every \(A \in \mathcal{I}\) and \(a, b \in \mathbb{R}\), the condition \(aA+b \in \mathcal{I}\) holds (e.g. \(\mathcal{I}\) is thw collection of Lebesgue null sets  \(\mathcal{N}\), or the collection of meager sets \(\mathcal{M}\)). Then \(\non(\mathcal{I}) = d(\mathbb{R}, \tau_\mathcal{I})\).
\end{example}
\begin{proof}
    Take \(A \not\in \mathcal{I}\) with \(|A| = \non(\mathcal{I})\). There exists a bounded interval \((a,b)\) such that \(A \cap (a,b) \not\in \mathcal{I}\) (otherwise, \(\bigcup\{[k, k+1) \cap A: k \in \mathbb{Z}\} = A \not\in \mathcal{I}\), while \([k, k+1) \cap A \in \mathcal{I}\), but \(\mathcal{I}\) is a \(\sigma\)-ideal). For any \(p < q\), let \(\phi_{p,q}: (a,b) \to (p,q)\) denote the natural increasing bijection. By the assumption on the ideal, \(\phi_{p,q}[A] \not\in \mathcal{I}\). Furthermore, let us define \(\tilde{A} = \bigcup_{p\in\mathbb{Q}}\bigcup_{q \in \mathbb{Q} \cap (p, \infty)} \phi_{p,q}[A]\). Notice that \(|\tilde{A}| = |A|\) and the set \(\tilde{A}\) is dense in \((\mathbb{R}, \tau_\mathcal{I})\): for any basic open set of the form \((p,q) \setminus B\) such that \(p,q \in \mathbb{Q}\), \(p < q\), and \(B \in \mathcal{I}\), we have \(((p,q) \setminus B) \cap \tilde{A} \supseteq ((p,q) \setminus B) \cap \phi_{p,q}[A] = \phi_{p,q}[A] \setminus B \not= \emptyset\), since \(\phi_{p,q}[A] \not\in \mathcal{I}\) and \(B \in \mathcal{I}\).
\end{proof}

These considerations suggest a natural question regarding the relationship between the cardinal invariants of an ideal space, those of the underlying space, and the invariants of the ideal used in its definition. In the remainder of this section, we explore this problem.

\begin{theorem}\label{covI}
Let \((X,\tau)\) be a \(T_1\) space, and let \(\mathcal{I}\) be an ideal on \(X\) to which all sets of cardinality not greater than \(t(X, \tau)\) belong. Then \(\cov(\mathcal{I}) \leq \max\{ w(X, \tau_\mathcal{I}), L(X, \tau)\}\).
\end{theorem}
\begin{proof}
    The set \(\{U \setminus A: U \in \tau, A \in \mathcal{I}\}\) is a base for the space \((X, \tau_\mathcal{I})\). It follows from the above lemma that there exists a base \(\mathcal{B} = \{U_\alpha \setminus A_\alpha: \alpha < w(X, \tau_ \mathcal{I})\}\). Let us define \(\mathcal{A} = \{A_\alpha: \alpha < w(X, \tau_\mathcal{I})\}\). We will show that there exists a set \(A\) of cardinality strictly less than \(t(X, \tau) + L(X, \tau)\) such that \(A \cup \bigcup \mathcal{A} = X\).
    
    Suppose for the sake of contradiction that \(X \setminus \bigcup \mathcal{A}\) has cardinality at least \(t(X, \tau) + L(X, \tau)\). Since \(|X \setminus \bigcup \mathcal{A}| \geq L(X, \tau)\), the set \(X \setminus \bigcup \mathcal{A}\) has an accumulation point \(x\). Furthermore, since \(|X \setminus \bigcup \mathcal{A}| \geq t(X, \tau)\), we can choose a set \(E\) such that \(x \not\in E\) and \(E \cap \bigcup \mathcal{A} = \emptyset\). The set \(X \setminus E\) is open in \((X, \tau_\mathcal{I})\). Since \(\mathcal{B}\) is a base, there exists \(S \subseteq w(X, \tau_\mathcal{I})\) such that \(X \setminus E = \bigcup \{U_\alpha \setminus A_\alpha: \alpha \in S\}\). Let \(U = \bigcup \{U_\alpha: \alpha \in S\}\). Notice that \(X \setminus E \subseteq U\).
    
    If \(E \cap U = \emptyset\), then \(U \subseteq X \setminus E\), hence \(U = E\). However, \(U \in \tau\), which means that \(X \setminus U = E\) is a closed set in \((X, \tau)\). This contradicts the fact that \(x \not\in E\).
    On the other hand, if \(x' \in U \cap E\), then we can find \(\alpha_0\) such that \(x' \in U_{\alpha_0}\). By definition, we have \(x' \not\in A_{\alpha_0}\), which leads to a contradiction since
    \[
    x' \in U_{\alpha_0} \setminus A_{\alpha_0} \subseteq \bigcup \{U_\alpha \setminus A_\alpha: \alpha \in S\} \subseteq X \setminus \{x'\}.
    \]

    We have shown that \(|X \setminus \bigcup \mathcal{A}| < t(X, \tau) + L(X, \tau) = t(X, \tau) \cdot L(X, \tau)\). Let \(\mathcal{A}'\) be a partition of \(X \setminus \bigcup \mathcal{A}\) into at most \(L(X, \tau)\) sets of cardinality not greater than \(t(X, \tau)\). Then \(\mathcal{P} = \mathcal{A} \cup \mathcal{A}' \subseteq \mathcal{I}\) and \(\mathcal{P}\) covers the set \(X\). It follows that
    \[
    \cov(\mathcal{I}) \leq |\mathcal{P}| \leq w(X, \tau_\mathcal{I}) + L(X, \tau).
    \]
\end{proof}

\begin{theorem}
Under the assumptions of theorem \ref{covI}:
\begin{enumerate}
    \item \(\Con(\cov(\mathcal{I}) > d(X, \tau_\mathcal{I}))\),
    \item \(\Con(\cov(\mathcal{I}) < w(X, \tau_\mathcal{I}))\).
\end{enumerate}
\end{theorem}
\begin{proof}
Let \((\mathbb{R}, \tau)\) denote the real line with the standard topology.
\begin{enumerate}
    \item Suppose that \(\mathfrak{c} = \aleph_2\). Let \(\mathcal{I}\) be the ideal of countable subsets of \(\mathbb{R}\). Furthermore, for any \(x, y \in \mathbb{R}\), let us define \(x \sim y \iff x-y \in \mathbb{Q}\). Notice that \(\sim\) is an equivalence relation. Take a set \(\mathcal{A} \subseteq \mathbb{R}/\!\sim\) of cardinality \(\aleph_1\).

    The set \(\bigcup \mathcal{A}\) is dense in \((\mathbb{R}, \tau_\mathcal{I})\). Indeed, if \(U \not= \emptyset\) is open in \((X, \tau)\), then \(|U \cap \bigcup \mathcal{A}| = \aleph_1\), so for any countable set \(I\), we have \((U \setminus I) \cap \bigcup \mathcal{A} \not= \emptyset\). We therefore have
    \[
    d(\mathbb{R}, \tau_\mathcal{I}) = \aleph_1 < \aleph_2 = \cov(\mathcal{I}).
    \]

    \item Suppose that \(\mathfrak{c} = \aleph_{\omega_1}\). Let \(\mathcal{I} = \{A \subseteq \mathbb{R}: |A| < \mathfrak{c}\}\). Then \(\cov(\mathcal{I}) = \aleph_1 < \mathfrak{c} = d(\mathbb{R}, \tau_\mathcal{I})\).
\end{enumerate}
\end{proof}

The inequality \(\cov(\mathcal{I}) < w(X, \tau_\mathcal{I})\) can be obtained without the assumption of the singularity of \(\mathfrak{c}\). It suffices to take the set of real numbers with the standard topology and the ideal \(\mathcal{I} = \{A \subseteq \mathbb{R}: \inf A > -\infty \wedge \sup A < \infty\}\). Then \(\cov(\mathcal{I}) = \aleph_0 < \mathfrak{c} = w(X, \tau_\mathcal{I})\).

\begin{theorem}
For any topological space \((X, \tau)\) and any ideal \(\mathcal{I}\) on \(X\), the inequality \(t(X, \tau_\mathcal{I}) \leq \max\{\cof(\mathcal{I}), \chi(X, \tau)\}\) holds.
\end{theorem}
\begin{proof}
Take any point \(x \in X\). We will show that the tightness at this point is at most \(\cof(\mathcal{I}) \cdot \chi(X, \tau)\). To this end, take \(Y \subseteq X\) such that \(x \in \cl^*(Y) = Y \cup Y^*\). We want to find \(Z \subseteq Y\) such that \(x \in \cl^*(Z)\) and \(|Z| \leq \cof(\mathcal{I}) \cdot \chi(X, \tau)\).

If \(x \in Y\), then \(Z = \{x\}\) suffices. Suppose that \(x \in Y^*\). This means that for every \(U \in \tau(x)\), we have \(U \cap Y \not\in \mathcal{I}\). For \(A \in \mathcal{I}\), we have \(U \cap Y \not\subseteq A\). By the axiom of choice, for any \(U \in \tau(x)\) and \(A \in \mathcal{I}\), let us choose \(x_{UA} \in (U \cap Y) \setminus A\). Define \(Z = \{x_{UA}: U \in \mathcal{B}, A \in \mathcal{A}\}\), where \(\mathcal{B}\) is a neighbourhood base at \(x\) of minimal cardinality, and \(\mathcal{A}\) is a cofinal family in \(\mathcal{I}\) of minimal cardinality.

It is easy to see that \(Z \subseteq Y\) and \(|Z| \leq |\mathcal{B}| \cdot |\mathcal{A}| = \chi(X, \tau) \cdot \cof(\mathcal{I})\). It remains to show that \(x \in \cl^*(Z)\), i.e., \(\forall U \in \tau(x) ~ U \cap Z \not\in \mathcal{I}\). Equivalently, \(\forall U \in \mathcal{B} ~ U \cap Z \not\in \mathcal{I}\), or yet another way: \(\forall U \in \mathcal{B} ~ \forall A \in \mathcal{I} ~ U \cap Z \not\subseteq A\). Take \(U \in \mathcal{B}\) and \(A \in \mathcal{I}\). There exists \(A' \in \mathcal{A}\) such that \(A \subseteq A'\). We have \(x_{UA'} \in (U \cap Z) \setminus A' \subseteq (U \cap Z) \setminus A\), which concludes the proof.
\end{proof}

\begin{corollary}
If \((X, \tau)\) is first-countable, then \(t(X, \tau_\mathcal{I}) \leq \cof(\mathcal{I})\).
\end{corollary}

\section{Continuous functions}

It is easy to see that 
\[
C((X, \tau), (X, \tau_\mathcal I)) \subseteq C((X, \tau), (X, \tau)) \subseteq C((X, \tau_\mathcal I), (X, \tau))
\]
for any topological space \((X, \tau)\) and any ideal \(\mathcal I\) on \(X\). In this section we show that under the assumption that the space \((X, \tau)\) is well-behaved it is possible to prove much more about the relation between these three sets.

\begin{theorem}
Let \((\mathbb{R}, \tau)\) denote the real line with the Euclidean topology. Furthermore, let \(\mathcal{I}\) be the ideal of countable sets. Then every continuous function \(f: (\mathbb{R}, \tau) \rightarrow (\mathbb{R}, \tau_\mathcal{I})\) is constant.
\end{theorem}
\begin{proof}
Suppose for the sake of contradiction that \(f\) is a continuous and non-constant function. There exist \(a<b\) such that \(f(a) \not= f(b)\). Let us define for \(x \in \mathbb{R}\)
\[
g(x)=
\begin{cases*}
f(a), & if $x \leq a$, \\
f(x), & if $x \in (a, b)$, \\
f(b), & if $x \geq b$. \\
\end{cases*}
\]
Let \(U\) be an open set in \((\mathbb{R}, \tau_\mathcal{I})\). Notice that
\[
g^{-1}(U)=
\begin{cases*}
f^{-1}(U),                      & if $f(a) \not\in U \wedge f(b) \not\in U$, \\
f^{-1}(U) \cup (-\infty, a), & if $f(a) \in U \wedge f(b) \not\in U$, \\
f^{-1}(U) \cup (b, \infty),  & if $f(a) \not\in U \wedge f(b) \in U$, \\
f^{-1}(U) \cup (-\infty, a) \cup (b, \infty),  & if $f(a) \in U \wedge f(b) \in U$.
\end{cases*}
\]
In each case, the preimage of \(U\) under the function \(g\) is open, hence \(g\) is a continuous function.

Take any \(x_0 \in \mathbb{R}\) such that \(f(a) \not= g(x_0) \not= f(b)\). Let \(X_0 = g^{-1}(x_0)\) and \(m_0 = \max X_0\). The definition of \(m_0\) is well-posed, because by the continuity of \(g\), \(X_0\) is a closed set. If \(x_n, X_n, m_n\) are defined, let us take any \(x_{n+1} \in (m_n, b)\) such that \(f(a) \not= g(x_{n+1}) \not= f(b)\). We also define \(X_{n+1} = g^{-1}(g(x_{n+1}))\) and \(m_{n+1} = \max X_{n+1}\).

The sequence \((m_n)\) is monotonic and bounded, and hence convergent to some \(m \in \mathbb{R}\). Since the set \(Y=\{f(x_n): n \in \omega\}\) is countable, it is also closed in the space \((\mathbb{R}, \tau_\mathcal{I})\). By the continuity of \(g\), \(g^{-1}(Y)\) is a closed set in the standard topology. For every natural number \(n\), on the one hand \(m_n \in g^{-1}(Y)\), and on the other hand \(m \not\in X_n\), so \(m \not\in g^{-1}(Y)\). This contradicts the closedness of \(g^{-1}(Y)\).
\end{proof}

The conclusion of the theorem stated above is quite powerful. In fact, \(\mathbb R\) is not special here: if we replace \(\mathbb R\) with an arbitrary Polish space, we obtain a similar result, though with a slightly weaker conclusion. First, we need the following lemma. Here by \(\rho\) we mean the standard Polish topology on the Cantor space \(2^\omega\).

\begin{lemma}
\label{cantorToPolish}
Let \((X, \tau)\) be any Polish space, and let \(\mathcal{I} = \{A \subseteq X: |A|\leq\aleph_0\}\). Every continuous function \(f: (2^\omega, \rho) \rightarrow (X, \tau_\mathcal{I})\) has a finite set of values.
\end{lemma}
\begin{proof}
Suppose for the sake of contradiction that \(f\) takes distinct values at each of the points \(x_0, x_1, \dots\). By compactness, \((x_n)_{n\in\omega}\) has a convergent subsequence, so we can assume without loss of generality that \((x_n)_{n\in\omega}\) converges to \(x\). The set \(A=\{f(x_n): n\in\omega \wedge f(x_n) \not= f(x)\}\) is countable, hence closed in \((X, \tau_\mathcal{I})\). However, \(\forall^\infty n (x_n \in f^{-1}(A))\). Since \(f\) is continuous, \(f^{-1}(A)\) is closed, and thus \(x = \lim_{n\to\infty} x_n \in f^{-1}(A)\). Contradiction.
\end{proof}

Notice that for any \(n \in \omega\), the function \(f_n: 2^\omega \to \mathbb R\) defined as \(f_n(x) = \min(\{k: x(k)=1\}\cup\{n\})\) is in \(C((2^\omega, \rho), (X, \tau))\) and its range has cardinality \(n+1\). This shows that the previous lemma is somehow optimal.

\begin{theorem}
\label{XtoXI}
    Every continuous function \(f \in C((X, \tau), (X, \tau_\mathcal{I}))\) takes countably many values for any Polish space \((X, \tau)\) and \(\mathcal{I} = \{A \subseteq X: |A|\leq\aleph_0\}\).
\end{theorem}
\begin{proof}
Suppose for the sake of contradiction that \(f\) takes uncountably many values. By the Cantor-Bendixson theorem, it follows without loss of generality that \(X\) is a perfect space (equipped with a complete metric \(d\)). Let \(A=\{x_\alpha: \alpha < \omega_1\}\) be such that \(\alpha \not= \beta \implies f(x_\alpha) \not=f(x_\beta)\). In a perfect Polish space, every uncountable set has a condensation point. Let \(x_\infty\) denote the condensation point of \(A\). Let us denote \(C_\emptyset = \cl(B(x_\infty, 1))\).

If for \(n \in \omega\) and a sequence \(s\) of length \(n\) of zeros, the set \(C_s\) is defined, we define \(C_{s \concat 0}\) and \(C_{s \concat 1}\) as follows. Choose \(\alpha < \omega_1\) such that \(x_\alpha \in C_s\). Next, let \(r = \frac{1}{3} \min\{d(x_\infty, x_\alpha), \frac{1}{2^n}\}\) and \(C_{s \concat 0} = \cl(B(x_\infty, r))\), \(C_{s \concat 1} = \cl(B(x_\alpha, r))\). Note that \(C_{s \concat 0}\) and \(C_{s \concat 1}\) are disjoint closed sets.

The set \(C_{s \concat 1}\) has a subset homeomorphic to the Cantor set, of which \(x_\alpha\) is an element. Let \(h: 2^\omega \to C_{s \concat 1}\) be this embedding. For \(\sigma \in 2^\omega\), let us denote \(C_{s \concat 1 \concat \sigma} = h(\sigma)\).

Finally, let \(C_{(0,0,\dots)}=x_\infty\). The set \(\{C_\sigma: \sigma \in 2^\omega\}\) is homeomorphic to \(2^\omega\) and the function \(f\) takes infinitely many values on this set.
\end{proof}

Again, it is clear that theorem \label{XtoXI} is optimal in the sense that one cannot substitute “countably many” with “finitely many”. To see this, define \(f: \omega^\omega \to \omega^\omega\) by letting, for every \(x \in \omega^\omega\), \(f(x) = \{(n, x(0)) : n \in \omega\}\). Then the range of \(f\) is infinite and \(f\) is continuous.

This brings to an end our investigation of continuous functions whose domain is Polish and codomain is idealized. Notice that while the preceding theorems deal specifically with the ideal of countable sets, the result extends readily: the ideal can be replaced by any \(\sigma\)-ideal that contains all singletons, which yields a finer codomain.

We now turn to the converse case, where the domain is an idealized space.

\begin{theorem}
    Let \((\mathbb{R}, \tau)\) denote the real line with the Euclidean topology. Furthermore, let \(\mathcal{I}\) be an ideal on \(\mathbb{R}\) such that \(\mathcal{I} \cap \tau = \{\emptyset\}\). Every continuous function \(f: (\mathbb{R}, \tau_\mathcal{I}) \to (\mathbb{R}, \tau)\) is also continuous as a function from the space \((\mathbb{R}, \tau)\) to the space \((\mathbb{R}, \tau)\).
\end{theorem}
\begin{proof}
Suppose that \(f\) is a counterexample. There exists, therefore, a point of discontinuity, i.e., some \(x\) and an interval \((a,b) \ni f(x)\) such that \(\forall \epsilon ~ f((x-\epsilon, x+\epsilon)) \not\subseteq (a,b)\). Let us denote \(\delta = \frac{1}{2}\min\{f(x)-a, b-f(x)\}\) and \(f^{-1}(a,b) = U \setminus A\), \(f^{-1}(f(x) - \delta, f(x) + \delta) = V \setminus B\), where \(U,V \in \tau\), \(A,B \in \mathcal{I}\), \(A \subseteq U\). Let us take a sequence \((x_n)_{n\in\omega}\) of elements of \(A\) such that \(f(x_n) \not\in (a,b)\) for every \(n \in \omega\) and \(\lim_{n \to \infty} x_n = x\).
Since \(V\) is an open neighborhood of \(x\) and the sequence \((x_n)_{n \in \omega}\) converges to \(x\), there exists \(n\) such that \(x_n \in V\). By continuity, \(f^{-1}((f(x_n)-\delta, f(x_n)+\delta)) = W \setminus C\) for some \(W \in \tau\) and \(C \in \mathcal{I}\). Notice that the intervals \((f(x)-\delta, f(x)+\delta)\) and \((f(x_n) - \delta, f(x_n)+\delta)\) are disjoint, hence \((V\setminus B) \cap (W \setminus C) = \emptyset\), which means \(V \cap W \subseteq B \cup C \in \mathcal{I}\). However, \(V\) and \(W\) are open neighborhoods of \(x_n\), so \(\emptyset \not= V \cap W \in \tau\). This contradicts the assumption that \(\mathcal{I} \cap \tau = \{\emptyset\}\).
\end{proof}

\end{document}